\documentclass[11pt,a4paper]{article}

\usepackage{amsmath,amsfonts,amssymb,latexsym,graphics,epsfig,url}
\usepackage{xcolor}
\usepackage{amsthm}
\usepackage[english]{babel}
\usepackage{mathdots}
\usepackage{graphicx}
\usepackage{soul}
\usepackage[utf8]{inputenc}
\usepackage{diagbox}
\usepackage{cancel}
\usepackage{float}              
\newtheorem{theorem}{Theorem}[section]
\newtheorem{proposition}[theorem]{Proposition}
\newtheorem{lemma}[theorem]{Lemma}

\newtheorem{example}[theorem]{Example}

\DeclareMathOperator{\Aut}{Aut}

\DeclareMathOperator{\diag}{diag}
\DeclareMathOperator{\dist}{dist}
\DeclareMathOperator{\Fix}{Fix}

\DeclareMathOperator{\spec}{sp}

\def\Z{\ns Z}

\def\f{\mbox{\boldmath $f$}}

\def\x{\mbox{\boldmath $x$}}

\def\vecv{\mbox{\boldmath $v$}}

\def\vecx{\mbox{\boldmath $x$}}

\def\vecphi{\mbox{\boldmath $\phi$}}
\def\vec0{\mbox{\boldmath $0$}}

\def\A{\mbox{\boldmath $A$}}
\def\B{\mbox{\boldmath $B$}}

\def\D{\mbox{\boldmath $D$}}

\def\L{\mbox{\boldmath $L$}}
\def\M{\mbox{\boldmath $M$}}

\def\U{\mbox{\boldmath $U$}}

\def\Z{\ns{Z}}

\def\U{\mbox{\boldmath $U$}}

\def\1{\mbox{\boldmath $1$}}

\def\Re{\mathbb R}
\def\Z{\mathbb Z}

\begin{document}

	\title{A general method to find the spectrum and eigenspaces of the $k$-token of a cycle, and
		2-token through continuous fractions
		\thanks{This research has been supported by
			AGAUR from the Catalan Government under project 2021SGR00434 and MICINN from the Spanish Government under project PID2020-115442RB-I00.
			The research of M. A. Fiol was also supported by a grant from the  Universitat Polit\`ecnica de Catalunya with references AGRUPS-2022 and AGRUPS-2023.}}
	
	\author{M. A. Reyes$^a$, C. Dalf\'o$^a$, M. A. Fiol$^b$, and A. Messegu\'e$^a$\\
		{\small $^a$Dept. de Matem\`atica, Universitat de Lleida, Lleida/Igualada, Catalonia}\\
		{\small {\tt \{monicaandrea.reyes,cristina.dalfo\}@udl.cat}, {\tt  arnau.messegue@upc.edu}}\\
		{\small $^{b}$Dept. de Matem\`atiques, Universitat Polit\`ecnica de Catalunya, Barcelona, Catalonia} \\
		{\small Barcelona Graduate School of Mathematics} \\
		{\small  Institut de Matem\`atiques de la UPC-BarcelonaTech (IMTech)}\\
		{\small {\tt miguel.angel.fiol@upc.edu}}}

	\date{}
	\maketitle
	
	\begin{abstract}
		The $k$-token graph $F_k(G)$ of a graph $G$ is the graph whose vertices are the $k$-subsets of vertices from $G$, two of which being adjacent whenever their symmetric difference is a pair of adjacent vertices in $G$.
		In this paper, we propose a general method to find the spectrum and eigenspaces of the $k$-token graph $F_k(C_n)$ of a cycle $C_n$. The method is based on the theory of lift graphs and the recently introduced theory of over-lifts. 
		In the case of $k=2$, we use continuous fractions to derive the spectrum and eigenspaces of the 2-token graph of $C_n$. 
	\end{abstract}
	
	
	
	\section{Introduction}
	\label{sec:-1}
	
\begin{figure}[t]
		\begin{center}
			\vskip -2cm
			\includegraphics[width=14cm]{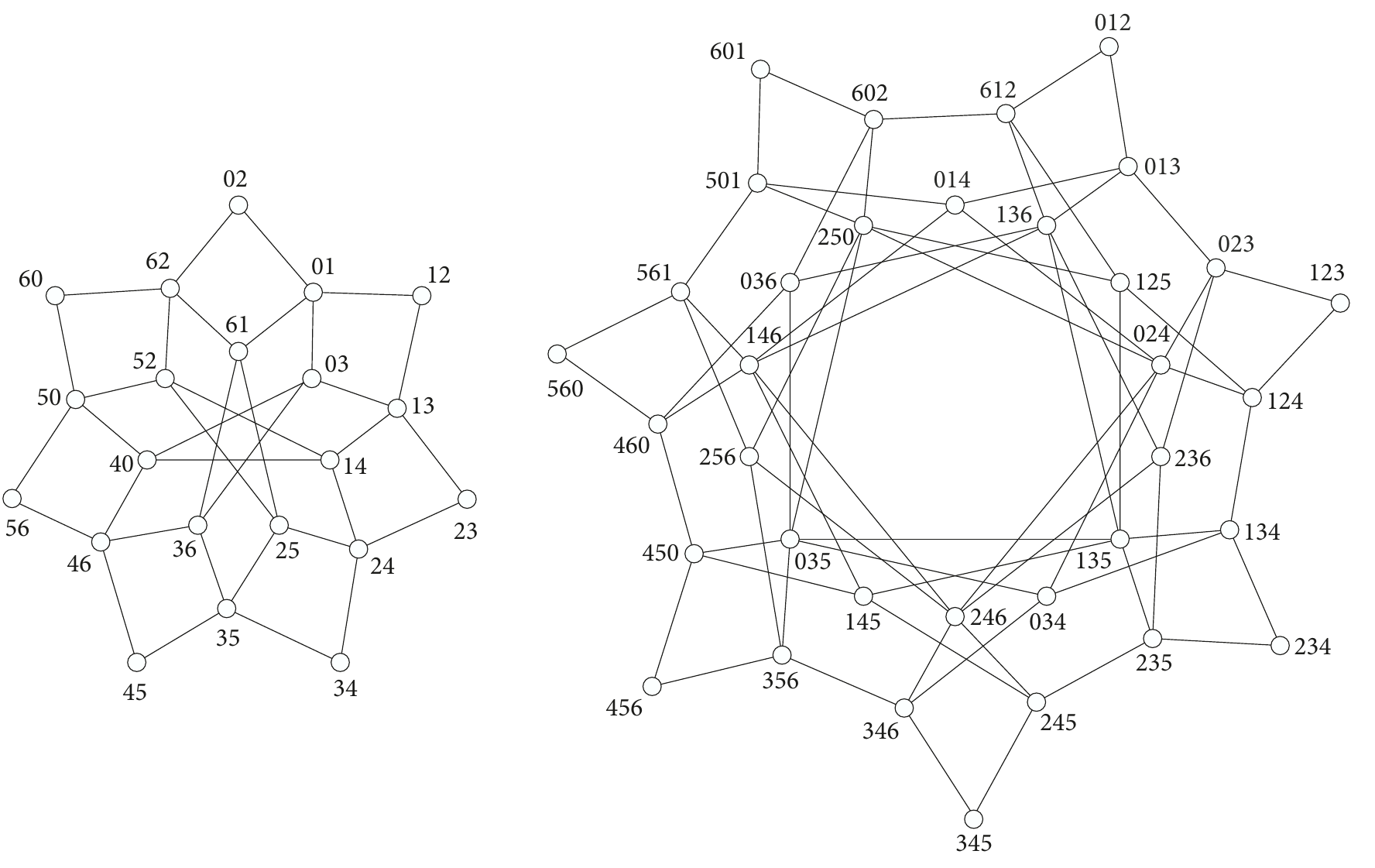}
			\caption{The 2-token and 3-token graphs of the cycle  $C_7$. }
			\label{F2(C7)+F3(C7)}
		\end{center}
	\end{figure}

	Let $G$ be a simple graph with vertex set $V(G)=\{1,2,\ldots,n\}$ and edge set $E(G)=\{(a,b):a,b\in V\}$. 
	For a given integer $k$ such that
	$1\leq k \leq n$, the {\em $k$-token graph} $F_k(G)$ of $G$ is the graph whose vertex set 
	$V(F_k(G))$ consists of the ${{n}\choose{k}}$
	$k$-subsets of vertices of $G$, and two vertices $A$ and $B$
	of $F_k(G)$ are adjacent whenever their symmetric difference $A \bigtriangleup B$ is a pair $\{a,b\}$ such that $a\in A$, $b\in B$, and $(a,b)\in E(G)$.
	For example, Figure \ref{F2(C7)+F3(C7)} shows the 2-token and 3-token graphs of the cycle $C_7$, where each vertex $\{a,b\}$ or $\{a,b,c\}$ is 
	represented by $ab$ or $abc$, respectively.
	
	The name `token graph' comes from an observation in
	Fabila-Monroy,  Flores-Pe\~{n}alo\-za,  Huemer,  Hurtado,  Urrutia, and  Wood \cite{ffhhuw12}, that vertices of $F_k(G)$ correspond to configurations of $k$ indistinguishable tokens placed at distinct vertices of $G$, where
	two configurations are adjacent whenever one configuration can be reached
	from the other by moving one token along an edge from its current position
	to an unoccupied vertex. 
	
	Note that if $k=1$, then $F_1(G)\cong G$; and if $G$ is the complete graph $K_n$, then $F_k(K_n)\cong J(n,k)$, where $J(n,k)$ denotes the Johnson graph~\cite{ffhhuw12}.
	
	
	Token graphs have some applications in physics. For instance, a
	relationship between token graphs and the exchange of Hamiltonian operators in quantum mechanics is given in Audenaert, Godsil, Royle, and Rudolph \cite{agrr07}.
	Moreover, token graphs are used for studying the isomorphism graph problem (see Fabila-Monroy and Trujillo-Negrete \cite{ft22}), and some error-correcting codes (see G\'omez Soto, Lea\~{n}os, R\'{i}os-Castro, and Rivera \cite{glrr18}).
	
	In this paper, we concentrate on the Laplacian spectrum and eigenspaces of the $k$-token $F_k(C_n)$ of a cycle $C_n$ for any value of $k$.
	However, our results can be applied to a universal adjacency matrix, that is, a
	linear combination with real coefficients of the adjacency matrix, the diagonal matrix of vertex degrees, the identity matrix, and the all-1 matrix. For more information, see Dalf\'o, Fiol, Pavl\'{\i}kov\'a, and \v{S}ir\'an
	\cite{df22b}. 
	Recall
	that the Laplacian matrix $\L(G)$ of a graph $G$ is $\L(G) =\D(G) -\A(G)$, where $A(G)$ is
	the adjacency matrix of $G$, and $\D(G)$ is the diagonal matrix whose diagonal entries are the vertex degrees of $G$.
	The Laplacian matrix is used because of its good properties (and
	applications) concerning token graphs. To mention just a few, the spectrum of $F_h(G)$ is contained in the spectrum of $F_h(G)$ for any $h\le k\le n/2$. This was proved by Dalf\'o, Duque, Fabila-Monroy, Fiol, Huemer, Trujillo-Negrete, and Zaragoza Mart\'{\i}nez \cite{ddffhtz21} by using a matrix analysis. Besides, as a consequence of the proof of Aldous' spectral gap conjecture, see 
	Caputo, Liggett, and Richthammer \cite{clr10} and Cesi \cite{c16}, the algebraic connectivity (that is, the smallest non-zero eigenvalue of the Laplacian) is the same for all $k$-token graph of a given graph.  Also, the authors of \cite{ddffhtz21} derived a close relationship between the spectra of the $k$-token graph of a graph $G$ and the $k$-token graph of its complement $\overline{G}$.
	
	
	The contents of this paper are the following.
	In Section \ref{resultats-coneguts}, we recall some theoretical results used in our study. Section 3 provides a general method to find the spectrum of the $k$-token graph $F_k(C_n)$ of a cycle $C_n$. This method is based on the theories of lift graphs and over-lifts described in Section \ref{resultats-coneguts}. In Section \ref{sec:cont.frac}, we use continuous fractions to derive the spectra and eigenspaces of $2$-token graphs of a cycle.
	Finally, in Section \ref{sec:charac-pol}, we compute the characteristic polynomial of $F_2(C_n)$.
	
	\section{Lift graphs and over-lifts}
	\label{resultats-coneguts}
	
	Let $\cal{G}$ be a group. An ({\em ordinary\/}) {\em voltage assignment} on the (di)graph $G=(V,E)$ (a graph or digraph)  is a mapping $\alpha: E\to \cal{G}$ with the property that $\alpha(a^-)=(\alpha(a^+))^{-1}$ for every pair of opposite arcs arc $a^+,a^-\in E$. Thus, a voltage assigns an element $g\in \cal{G}$ to each arc of the (di)graph, so that two mutually reverse arcs $a^+$ and $a^{-}$, forming an undirected edge, receive mutually inverse elements $g$ and $g^{-1}$. Then, the (di)graph $G$ and the voltage assignment $\alpha$ determine a new (di)graph $G^{\alpha}$, called the {\em lift} of $G$, which is defined as follows. The vertex and arc sets of the lift are simply the Cartesian products $V^{\alpha}=V\times \cal{G}$ and $E^{\alpha}=E\times \cal{G}$, respectively. Moreover, for every arc $a\in E$ from a vertex $u$ to a vertex $v$ for any $u,v\in V$ (possibly, $u=v$) in $G$, and for every element $g\in \cal{G}$, there is an arc $(a,g)\in E^{\alpha}$ from the vertex $(u,g)\in V^{\alpha}$ to the vertex $(v,g\alpha(a))\in V^{\alpha}$.
	The interest of this construction is that, from the base graph $G$ and the voltages, we can deduce some properties of the lift graph $G^{\alpha}$. This is usually done through a matrix associated with $G$ that, in the case when the group $\cal{G}$ is cyclic, every entry of such a matrix is  a polynomial in $z$ with integer coefficients, and so it is called 
	the associated `polynomial matrix' $\B(z)$. 
	Then, the whole (adjacency or Laplacian) spectrum and eigenspaces of $G^{\alpha}$ can be retrieved from $\B(z)$. 
	More precisely,  Dalf\'o, Fiol, Miller, Ryan, and \v{S}ir\'a\v{n} \cite{dfmrs17} proved the following result.
	\begin{theorem}[\cite{dfmrs17}]
		\label{th:sp-lifts}
		Let $R(n)$ be the set of $n$-th roots of unity, and consider the base graph $G=(V,E)$ with voltage assignment $\alpha$ with the cyclic group  ${\cal G}=\Z_n$. If $\vecx=(x_u)_{u\in V}$ is an eigenvector of $\B(z)$ with eigenvalue $\lambda$, then the vector $\vecphi=(\phi_{(u,j)})_{(u,j)\in V^{\alpha}}$ with 
		components $\phi_{(u,j)}=z^j x_u$ is an eigenvector of the lift graph $G^{\alpha}$ corresponding to the eigenvalue $\lambda$.
		Moreover, all the eigenvalues (including multiplicities) of $G^{\alpha}$ are obtained:
		$$
		\spec G^{\alpha} = \bigcup_{z\in R(n)}\spec(\B(z)).
		$$
	\end{theorem}
	For more information on lift graphs and digraphs, see 
	Dalf\'o, Fiol, Pavl\'ikov\'a, and  \v{S}ir\'an \cite{df22b}, 
	and Dalf\'o, Fiol, and \v{S}ir\'a\v{n} \cite{dfs19}.
	
	Notice that, generally, a given graph $H$ can be constructed as a lift graph, say $G^{\alpha}$, of a `small' base graph $G$ when the automorphism group $\Aut H$ is `big enough'. In other words, $H$ must possess some symmetries.
	However, sometimes this is not the case, as we will see when $H$ is the $k$-token of an $n$-cycle with $k$ a non-trivial divisor of $n$.
	To overcome this drawback, we introduce a new technique, already implicitly introduced in \cite{rdfm23}, that consists of `forgetting' the base graph and constructing a proper polynomial matrix $\B$ directly. Our aim is that such a matrix must also contain information about the spectrum of $H$. In some cases, $\B$ yields a small number of eigenvalues not present in the spectrum of $H$ (which correspond to vectors $\vecphi$ in Theorem \ref{th:sp-lifts} that are not eigenvectors of $H$), but our method identify these eigenvalues. This is the reason why we called this approach the method of `over-lifts'.
	
	\section{A general method for computing the spectrum of $F_k(C_n)$}
	In this section, we present a method to derive the spectrum of the $k$-token graph of an $n$-cycle. The method can be applied for any values of $n$ and $k$, and it is based on the theory of lift graphs and over-lifts.
	
	For given $n$ and $k\le n/2$, let $[0,n-1]=\{0,1,\ldots,n-1\}$ be the vertex set of the cycle $C_n$, and $V=\left[^n_k\right]$ the vertex set $V$ of $F_k(C_n)$, with $|V|={n\choose k}$. The following lemma gives the order of the base graph or polynomial matrix of $F_k(C_n)$ seen as a lift graph or over-lift.
	
	\begin{lemma}
		\label{lem:Aut}
		The automorphism group of the token graph $F_k(C_n)=(V,E)$ has a transitive subgroup ${\cal G}$ isomorphic to the cyclic group $\Z_n$. For every $r=0,1,\ldots, n-1$, let $d(r)=\gcd(n,r)$ and $o(r)=\frac{n}{d(r)}$.
		Then, the number $T(n,k)$ of orbits of ${\cal G}$ acting on $V=\left[^n_k\right]$ is
		\begin{equation}
			T(n,k)=|V/{\cal G}|=\frac{1}{n}\sum_{o(r)|k} {d(r)\choose k/o(r)}.
			\label{necklaces}
		\end{equation}
	\end{lemma}
	\begin{proof}
		Every vertex $A$ of $F_k(C_n)$ corresponds to $k$ tokens placed in different vertices of $C_n$, $A=\{a_1,\ldots,a_k\}$. Thus, for any integer $\gamma$, the mapping 
		$$
		A=\{a_1,\ldots,a_k\}\mapsto A+\gamma=\{a_1+\gamma,\ldots,a_k+\gamma\},
		$$ all arithmetic modulo $n$, is clearly an automorphism of $F_2(C_n)$. (In particular, it is known that
		$\Aut F_2(C_n)=\Aut C_n$, see Ibarra and Rivera \cite{ir22}.) 
		
		To prove \eqref{necklaces}, notice that, in fact, we want to compute
		the number of distinct necklaces with $k$ black beads (representing the position of the tokens) and $n-k$ white beads (vertices with no token). (Two necklaces are equivalent if a given rotation of the other can obtain one.) Although, for every given $k$, a formula can be obtained from P\'olya's enumeration theorem (see below), we derive \eqref{necklaces} for completeness by using Burnside's lemma. Let ${\cal G}$ be a finite group that acts on a set $V$. For each $g\in {\cal G}$, let 
		$\Fix(g)=\{x\in V:g\ast x=x\}$ denote the set of elements in $V$ that are fixed by $g$. Burnside's lemma gives the following formula for the number of orbits, denoted $|V/{\cal G}|
		$:
		\begin{equation}
			|V/{\cal G}|=\frac{1}{n}\sum_{g\in{\cal G}} |\Fix(g)|.
			\label{Burnside}
		\end{equation}
		Then, looking at ${\cal G}=\Z_n$ as a group of permutations generated by $g=(123\ldots n)$ (cyclic notation), we notice that the element $g^r$ consists of $d(r)=\gcd(n,r)$ cycles of length $o(r)=n/d(r)$. Moreover, for a vertex in $V=\left[^n_k\right]$ to be fixed, the $k$ black beads must totally `fulfill' some of the (equal length) cycles so that $o(r)|k$. Since there are $d(r)$ cycles, the number of ways of doing so is ${d(r)\choose k/o(r)}$, which corresponds to the claimed value of $\Fix(g^r)$. For instance, if $n=8$, $k=4$, and $g^4=(15)(26)(37)(48)$, the number of ways $4$ beads can fulfill 2 cycles of length 2 is ${4\choose 2}=6$. This completes the proof.
	\end{proof}
	
	In Table \ref{taula:necklaces}, we listed some of the values obtained from \eqref{necklaces}, together with the corresponding reference in the On-Line Encyclopedia of Integer Sequences \cite{oeis} for the sequences with $k=2,\ldots,7$.
	As commented above, the sequence $(T(n,k): n\ge 1)$ can be obtained by using P\'olya enumeration, and corresponds to the  $k$-th column in the example of the sequence A047996 in \cite{oeis}, where we can find the following alternative formula for \eqref{necklaces}, 
	\begin{equation}
		T(n, k) = \frac{1}{n}\sum_{d|\gcd(n,k)} \phi(d){n/d\choose k/d},
		\label{necklaces2}
	\end{equation}
	where $\phi$ is the Euler function. (This was first proved by Gilbert and Riordan in \cite{gr61}.)
	This gives, for instance, the sequence for $T(n,8)$: 
	$1,1,5,15,43,99,217,429,810,1430,\ldots$ 
	In particular, if $k=p$, a prime, Hadjicostas 
	proved a generalization of a conjecture made by Sloane and Lang for $p=7$, and showed that $T(n,p)=\lceil{n\choose p}/n\rceil$, see the comments for the sequence A032192 in \cite{oeis}.
	Notice also that, if $\gcd(n,k)=1$, all the orbits have $n$ vertices and, hence, $|V/{\cal G}|={n\choose k}/n$. These correspond to `aperiodic' necklaces (not consisting of a repeated subsequence), and their number is given by the so-called `Moreau's necklace-counting function' (introduced by Moreau in 1872),
	\begin{equation}
		M(n, k) = \frac{1}{n}\sum_{d|\gcd(n,k)} \mu(d){n/d\choose k/d},
		\label{necklaces2-moreau}
	\end{equation}
	where $\mu$ is the classic M\"obius function. 
	In the case of periodic necklaces (with the period being the smallest length of the repeated subsequence), there are orbits with less than $n$ vertices. For instance, when $n=8$ and $k=4$, there are eight orbits with $n=8$ vertices, one orbit with $n/2=4$ vertices, and one orbit with $n/4=2$ vertices. 
	An efficient algorithm for generating fixed-density (that is, with a fixed number of zeros) $k$-ary necklaces or aperiodic necklaces was proposed by Ruskey and Savada \cite{rs99}.
	
	\begin{table}[t!]
		\begin{center}
			\begin{tabular}{|c|c|cccccccccc| }
				\hline
				$k\backslash n$ & OEIS & 3  & 4  & 5 & 6 & 7 & 8 & 9 & 10 & 11 & 12   \\
				\hline\hline
				2 & A004526 & 1 & 2 & 2 &  3 & 3 & 4 & 4 & 5 & 5 & 6   \\
				\hline
				3 & A007997 & 1 & 1 & 2 &  4 & 5 & 7 & 10 & 12 & 15 &  19 \\
				\hline
				4 & A008610 &  & 1  &  1 & 3 & 5 & 10 & 14 & 22 & 30 &  43 \\
				\hline
				5 & A008646 &   &    &   1 & 1 & 3 & 7 & 14 & 26 & 42 & 66  \\
				\hline
				6 & A032191 &     &   &   &  1 & 1 & 4 & 10 & 22 & 42 & 80 \\
				\hline
				7  & A032192 &   &   &   &   & 1 & 1 & 4 & 12 & 30 & 66\\
				\hline
			\end{tabular}
		\end{center}
		\caption{The number of orbits of the group $\Z_n$ acting on the $k$-token graph $F_2(C_n)$ or, alternatively, the number of necklaces with $k$ black beads and $n-k$ white beads, see the corresponding sequence in \cite{oeis}. }
		\label{taula:necklaces}
	\end{table}
	
\begin{figure}[t]
		\begin{center}
			\includegraphics[width=5.2cm,angle=0]{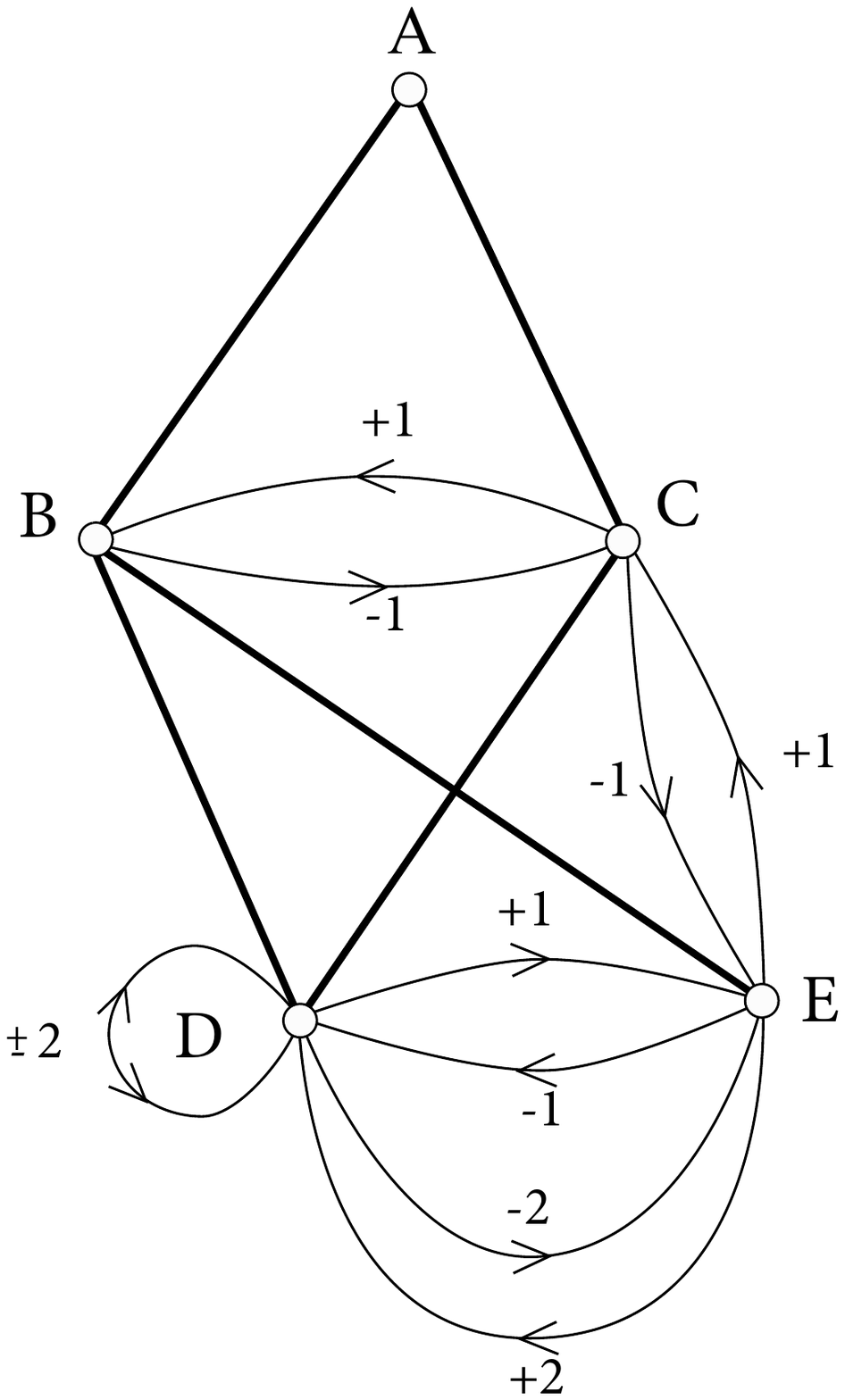}
			\caption{The base graph $G$ for the 3-token graph of $C_7$. }
			\label{base-graph}
		\end{center}
	\end{figure}

	\begin{theorem}
		\label{th:overlift}
		For every $n$ and $k\le n/2$, the spectrum of the token graph $F_k(C_n)$ can be obtained either from a lift graph of a base graph on the cyclic group $\Z_n$ or from an 
		over-lift polynomial matrix on the same group.
	\end{theorem}
	\begin{proof}
		Since, in both cases, the key information is given by the (Laplacian) polynomial matrix $\B$, with order $\nu$ given by Lemma \ref{lem:Aut}, it is enough to show how to obtain it.
		Then, the method goes along the following steps:
		\begin{enumerate}
			\item
			For each of the $\nu$ orbits, choose a representative vertex of $F_k(C_n)$, so constituting the set  ${\cal R}=\{A,B,C,\ldots\}$.
			\item
			Construct a $\nu\times\nu$, with $\nu=|{\cal R}|$, matrix $\B(z)$, indexed by the vertices $A,B,C,\ldots$, with diagonal entries the degrees $\deg(A),\deg(B),\ldots$ in $F_{k}(C_n)$.
			\item
			Now, suppose that $A$ is adjacent to $A'$ in $F_{k}(C_n)$. Then, by Lemma \ref{lem:Aut}, there exists some $B\in {\cal R}$ and $r\in \Z_n$ such that $A'=B+r$ (where, if $B=\{b_1,\dots,b_k\}$, then $B+r=\{b_1+r,\dots,b_k+r\}$). Then the polynomial entry $\B(z)_{AB}$ has the term $-z^{r}$ (in the putative base graph, the arc $A\rightarrow B$ would have voltage $+r$).
			\item
			Repeat the previous step for every vertex of ${\cal R}$ and its adjacent vertices in $F_k(C_n)$.
		\end{enumerate}
		In what follows, let ${\cal R}=\{A_1,\ldots,A_{\nu}\}$ be the set of indexes of $\B$ (or vertex set of the base graph). Let $\f= (f_1,\ldots,f_{\nu})$ be a $\lambda$-eigenvector of $\B$.
		Let $\L$ be the Laplacian matrix of $F_k(C_n)$. Then, for every  vertex $X$ of $F_k(C_n)$, we know that there exists $A_i\in{\cal R}$ and $j\in \Z_n$ such that $X=A_i+j$. Now, we claim that the vector $\vecv$ with ${n \choose k}$ components of the form
		\begin{equation}
			\label{vector-y}
			v_{X}=f_{i}z^j\qquad X\in {\cal R},\  j=0,\ldots,n-1,
		\end{equation}
		where $z=\zeta^r=e^{ir\frac{2\pi}{n}}$, is a $\lambda$-eigenvector of $\L$, provided that the following condition holds:
		\begin{itemize}
			\item[$(*)$]
			If $A_{i_1},\ldots,A_{i_t}$ are the periodic elements of ${\cal R}$, with respective periods $\pi(A_{i_1}),\ldots,\pi(A_{i_t})$, then, for every $j=1,\ldots,r$, either $f_{i_j}=0$, or $o(r)=\frac{n}{\gcd(n,r)}$ divides $\pi(A_{i_j})$.
		\end{itemize}
		To prove the claim, assume that, in $F_k(C_n)$, the vertex $A_i\in {\cal R}$, with degree $\delta$, is adjacent to the vertices $A_{i_1}+j_1,\ldots,A_{i_d}+j_d$. Then, the $i$-th equality of $\B\f=\lambda\f$ reads
		$$
		\delta f_i-\sum_{h=1}^{\delta} f_{i_h} z^{j_h}=\lambda f_i, 
		$$
		which corresponds to the $A_i$-th equation in $\L\vecv=\lambda\vecv$.
		Now, multiplying all terms by $z^p=\zeta^{rp}=e^{ipr\frac{2\pi}{n}}$, we get
		\begin{equation}
			\label{Lv=}
			\delta f_iz^p-\sum_{h=1}^{\delta} f_{i_h} z^{j_h+p}=\lambda f_i z^p,\qquad p=0,1,\ldots, n-1,
		\end{equation}
		which includes all equalities in $\L\vecv=\lambda\vecv$ corresponding to the vertices that are in the same orbit as $A_i$. Thus, in general, only the first $o(r)$ equations in \eqref{Lv=} matter (the following ones are repetitions). However,  since the above reasoning must hold for every $A_i\in {\cal R}$, to be consequent with   \eqref{Lv=}, we need to impose the conditions $(*)$ on $\f$ or $o(r)$.
		\\
		This concludes the proof.
	\end{proof}
	To clarify the procedure, let us see two examples. For simplicity, we represent each element $A=\{a_1,\ldots, a_k\}\in {\cal R}$ as a sequence $a_1\ldots a_k$, where the order of the digits is irrelevant. In the first example, the graph can be seen as a lift graph; hence, we equivalently construct the base graph with its voltages.
	In the second example, such a (`small enough') base graph does not exist, and then, we use the technique of over-lifts deriving the polynomial matrix directly.
	\begin{example}
		Let us first  consider the case of $F_3(C_7)$ shown on
		the right-hand side of Figure \ref{F2(C7)+F3(C7)}.
		Since $3\nmid 7$, the theorem guarantees that  $F_3(C_7)$ is the lift of a base graph $G$ on $\nu= {7\choose 3}/7=5$ vertices. Moreover, each of these vertices, say $A, B, C, D, E$, is a representative of the five orbits. For instance, 
		take $A=123$, $B=612$, $C=013$, $D=136$, and $E=125$.
		Then, the base graph $G$ has the following arcs and voltages:
		\begin{itemize}
			\item[$(A)$] 
			Vertex $A=012$ is adjacent to $612=B$ and $013=C$. Thus, the arcs $A\rightarrow B$ and $A\rightarrow C$
			have voltage $0$.
			\item[$(B)$] 
			Vertex $B=612$ is adjacent to $012=A$, $512=E$, $613=D$, and $602+1=013=C$ (the sum applies to every digit). Thus, the arcs $B\rightarrow A$,  $B\rightarrow E$, and $B\rightarrow D$ are arcs with voltage 0, whereas the arc $B\rightarrow C$ has voltage $-1$.
			\item[$(C)$] 
			Vertex $C=013$ is adjacent to $613=D$, $012=A$, $023-1=612=B$ and $014+1=125=E$. Therefore, the arcs $C\rightarrow D$ and  $C\rightarrow A$ have voltage $0$,  $C\rightarrow B$ has voltage $+1$, and  $C\rightarrow E$ has voltage $-1$.
			\item[$(D)$] 
			Vertex $D=613$ is adjacent to $013=C$, $612=B$, $513-2=361=D$, $623-1=512=E$, $603+2=125=E$, $614+2=136=D$. Therefore, the arcs $D\rightarrow C$ and  $D\rightarrow B$ have voltage $0$,  the two arcs $D\rightarrow D$ (loops) have voltages $\pm 2$, and the two arcs $D\rightarrow E$ have voltages $+1$ and $-2$.
			\item[$(E)$] 
			Vertex $E=125$ is adjacent to $126=B$, $025+1=613=D$, $135-2=613=D$, $124-1=013=C$. Therefore, the arc $E\rightarrow B$ has voltage $0$,  the two arcs $E\rightarrow D$  have voltages $-1$ and $+2$, and the two arc $E\rightarrow C$ has voltage $+1$.
		\end{itemize}
		The obtained base graph $G$ and its voltages are shown in Figure \ref{base-graph}.
	\end{example}
	Then, the polynomial matrix of $G$ is 
	$$
	\B(z)=\left(
	\begin{array}{ccccc}
		2 & -1 & -1 & 0 & 0\\
		-1 & 4 & -z^{-1} & -1 & -1\\
		-1 & -z & 4 & -1 & -z^{-1}\\
		0 & -1 & -1 & 6 - z^2 - z^{-2} & -z - z^{-2}\\
		0 & -1 & -z & -z^{-1} - z^2 & 4
	\end{array}
	\right).
	$$
	Since $F_3(C_7)$ is a lift graph, by Theorem \ref{th:sp-lifts},  its 
	${7\choose 3}=35$ eigenvalues can be obtained from $\B(z)$, see Table \ref{taula:F3(C7)}.
	\begin{table}[ht!]
		\begin{center}
			\begin{tabular}{|c|ccccc| }
				\hline
				$\zeta=e^{i\frac{2\pi}{7}}$, $z=\zeta^r$ & $\lambda_{r,1}$  & $\lambda_{r,2}$  & $\lambda_{r,3}$   & $\lambda_{r,4}$  & $\lambda_{r,5}$  \\
				\hline\hline
				$\spec(\B(\zeta^0))$ & \bf 0 & \bf 2.0 &  5.0 & 5.0 & \bf 6.0   \\
				\hline
				$\spec(\B(\zeta^1))=\spec(\B(\zeta^6))$ & \bf 0.7530   & 2.91929 & \bf 3.9363  &  5.7238 & \bf 7.1125 \\
				\hline
				$\spec(\B(\zeta^2))=\spec(\B(\zeta^5))$ &  \bf 1.1633 & \bf 2.4450 & 3.8385 & \bf 5.1446 &  9.2103 \\
				\hline
				$\spec(\B(\zeta^3))=\spec(\B(\zeta^4))$  & 1.2696 & \bf 1.9019 & \bf 3.8019  &  \bf 4.7411 & 7.0383 \\
				\hline
			\end{tabular}
		\end{center}
		\caption{All the eigenvalues of the matrices $\B(\zeta^r)$, which yield the eigenvalues of the 3-token graph $F_3(C_7)$. The values in boldface correspond to the eigenvalues of $C_7$ and $F_2(C_7)$.}
		\label{taula:F3(C7)}
	\end{table}
	\\
	\begin{table}[t]
		\begin{center}
			\begin{tabular}{|c|cccc| }
				\hline
				$\zeta=e^{i\frac{2\pi}{6}}$, $z=\zeta^r$ & $\lambda_{r,1}$  & $\lambda_{r,2}$  & $\lambda_{r,3}$   & $\lambda_{r,4}$  \\
				\hline\hline
				$\spec(\B(\zeta^0))$ & \bf \bf 0 &  \bf  2.7639 &  6.0 & \bf 7.2361  \\
				\hline
				$\spec(\B(\zeta^1))=\spec(\B(\zeta^5))$ & \bf 1.0   & 4.0 & \bf 5.0  &  $6.0^*$ \\
				\hline
				$\spec(\B(\zeta^2))=\spec(\B(\zeta^4))$ &  \bf 1.4384 & \bf 3.0 & $6.0^*$ &  \bf 5.5616 \\
				\hline
				$\spec(\B(\zeta^3))$  & 1.3944 & \bf 2.0 &  \bf 4.0 &  8.6056 \\
				\hline
			\end{tabular}
		\end{center}
		\caption{All the eigenvalues of the matrices $\B(\zeta^r)$, which yield the eigenvalues of the 3-token graph $F_3(C_6)$ plus four $6$'s (those marked with `*'). The values in boldface correspond to the eigenvalues of $C_6$ and $F_2(C_6)$.}
		\label{taula:C6}
	\end{table}
	\begin{figure}[t]
		\begin{center}
			\includegraphics[width=7cm]{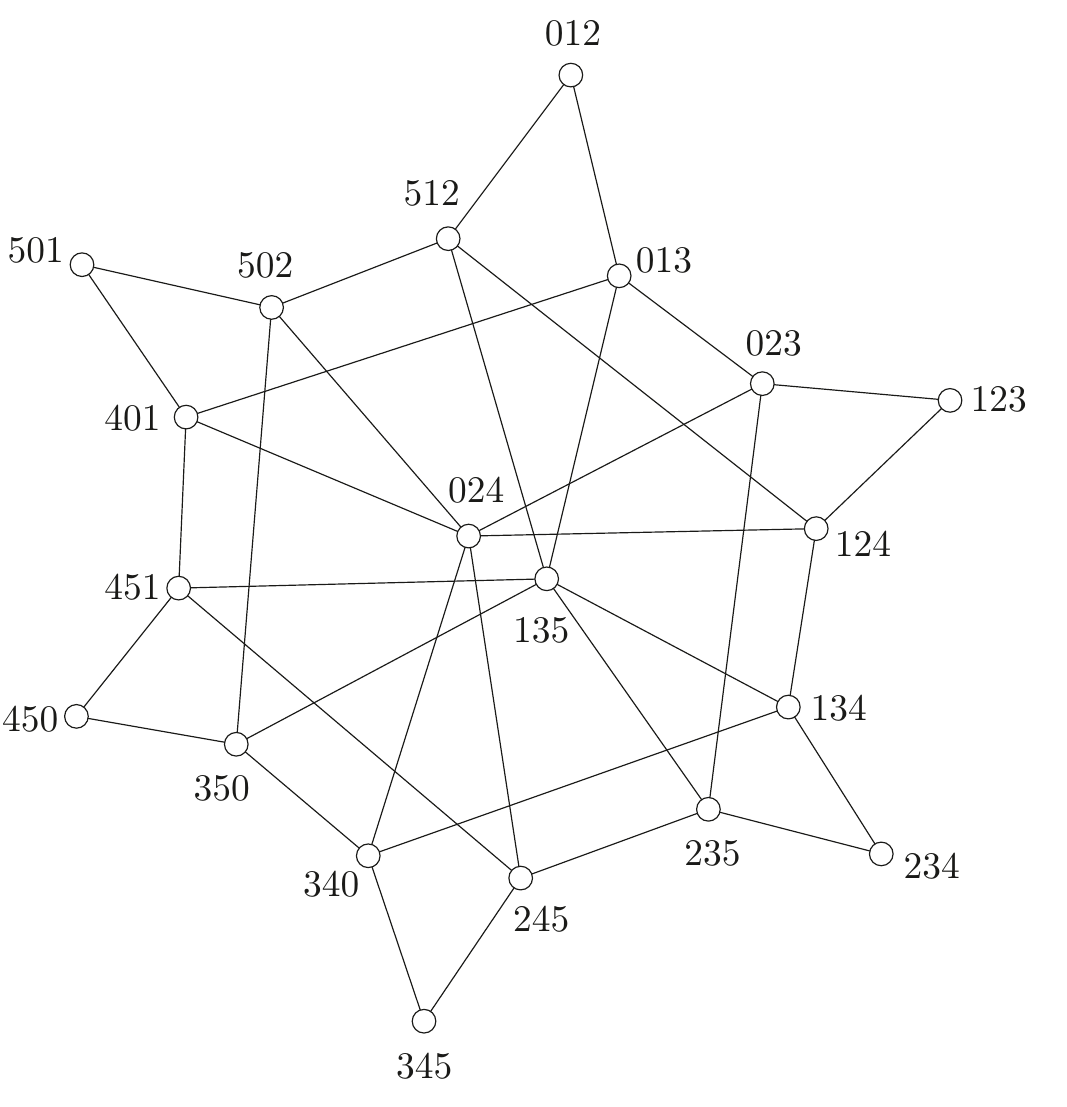}
			\caption{The 3-token graph of $C_6$.}
			\label{F3(C6)}
		\end{center}
	\end{figure}
	\begin{example}
		\label{ex:F3(C6)}
		Consider now the case of the token graph $F_3(C_6)$ shown in Figure \ref{F3(C6)}.
		Since $3|6$, we have 3 orbits with 6 vertices and one orbit with 2 vertices. Then, $F_3(C_6)$ can be obtained as an over-lift with polynomial $\nu\times \nu$ matrix of size $\nu=3+1$. As representatives of these orbits, we can take, for instance, 
		$A=012$, $B=013$, $C=014$,  and $D=024$.
		Then, we reason as follows:
		\begin{itemize}
			\item[$(A)$] 
			Vertex $A=012$ is adjacent to $512=C+1$ and $013=B$. Thus, $\B(z)_{AC}=-z$ and $\B(z)_{AB}=-1$.
			\item[$(B)$] 
			Vertex $B=013$ is adjacent to $513=D+1$, $023=C+2$, $012=A$, and $014=C$. Thus, $\B(z)_{BD}=-z$,   $\B(z)_{BC}=-1-z^2$, and  $\B(z)_{BA}=-1$.
			\item[$(C)$] 
			Vertex $C=014$ is adjacent to $514=B-2$, $024=D$, $013=B$, and $015=A-1$. Therefore,  $\B(z)_{CB}=-1-z^{-2}$,
			$\B(z)_{CD}=-1(\equiv -3,-5)$, and $\B(z)_{CA}=-z^{-1}$.
			\item[$(D)$] 
			Vertex $D=024$ is adjacent to $124=B+1$, $245=C+2$, $014=C$, $034=B+3$, $502=B+5$, $023=C+4$. Therefore, $\B(z)_{DB}=-z-z^3-z^{5}$ and $\B(z)_{DC}=-1-z^{2}-z^4$.
		\end{itemize}
		All the other off-diagonal entries of $\B(z)$ are zero, whereas the diagonal entries are $\deg(A)=2$, $\deg(B)=\deg(C)=4$, and $\deg(D)=6$.
		Summarizing, the polynomial matrix of $F_3(C_6)$ is:
		$$
		\B(z)=\left(
		\begin{array}{cccc}
			2 & -1 & -z & 0 \\
			-1 & 4 & -1-z^{2} & -z\\
			-z^{-1} & -1-z^{-2} & 4 & -1 \\
			0 & -z-z^3-z^5 & -1 - z^2 - z^{4} & 6\\
		\end{array}
		\right).
		$$
		Since $F_3(C_6)$ is an over-lift, its eigenvalues can also be obtained from $\B(z)$. However, according to the proof of Theorem \ref{th:overlift}, the spectrum of such a matrix contains some `spurious' 
		eigenvalues, not in the spectrum of  $F_3(C_6)$. As shown in Table  \ref{taula:C6}, such additional eigenvalues are four 6's. Following the reasoning in the proof of Theorem \ref{th:overlift}, the reason is that,
		for $r=1,5$, the $\lambda$-eigenvectors of $\B(\zeta^r)$ for $\lambda=6$
		are $(0.1\pm 0.173i, -0.2\mp 0.376i, -0.4,1)^{\top}$, and for $r=2,4$
		such $\lambda$-eigenvectors are  $(0.5\pm 0.865i, 1\mp 1.73i, -2,1)^{\top}$. Thus, since the last component is $f_4(=f_D)=1(\neq 0)$ and $\pi(D)=2$, we have that  $o(1)=o(5)=6\nmid \pi(D)$ and $o(2)=o(4)=3\nmid \pi(D)$. Consequently, none of the above $6$-eigenvectors  yields an eigenvector of (the Laplacian matrix of) $F_3(C_6)$. In contrast, for $r=0$ ($z=1$), $\B(1)$ has also an eigenvalue 6 with the corresponding eigenvector $(0,-1,1,0)^{\top}$. Then, since $f_4=0$, such an eigenvector also gives an eigenvector of $F_3(C_6)$, and the eigenvalue $6$ also belongs to its spectrum.
	\end{example}
	
	\begin{example}
		\label{ex:F4(C8)}
		Consider the case of the token graph $F_4(C_8)$.
		In this case, some vertices (seen as necklaces) are periodic because there are orbits with less than $n=8$ vertices. More precisely, there are eight orbits with $n=8$ vertices, one orbit with $n/2=4$ vertices, and one orbit with $n/4=2$ vertices. 
		Let $A=0123$, $B=1234$, $C=0125$, $D=0126$, $E=0134$, $F=0135$, $G=0136$, $H=0145$, $I=0146$,  and $J=0246$ be the vertices that represent each of these orbits. Then, the spectrum of  
		$F_4(C_8)$ can be obtained as an over-lift with polynomial $\nu\times \nu$ matrix of size $\nu=10$, where the polynomial matrix of $F_4(C_8)$ is
		$$
		\B(z)=
		\scriptsize
		\left(
		\arraycolsep=1pt
		\begin{array}{cccccccccc}
			2 & -1 & 0 & -z & 0 & 0 & 0 & 0 & 0 & 0 \\
			-1 & 4 & -1 & 0 & -1 & 0 & -z & 0 & 0 & 0 \\
			0 & -1 & 4 & -1 & 0 & -1 & 0 & 0 & -z & 0 \\
			-z^{-1} & 0 & -1 & 4 & -z^{-2} & 0 & -1 & 0 & 0 & 0 \\
			0 & -1 & 0 & -z^2 & 4 & -1 & 0 & 0 & -z^3 & 0 \\
			0 & 0 & -1 & 0 & -1 & 6 & - 1-z^2 & -1 & 0 & -z \\
			0 & -z^{-1} & 0 & -1 & 0 & -1 - z^{-2} & 6 & 0 & -z^2 - 1 & 0 \\
			0 & 0 & 0 & 0 & 0 & -z^4 - 1 & 0 & 4 & -z^4 - 1 & 0 \\
			0 & 0 & -z^{-1} & 0 & -z^{-3} & 0 & -1 - z^{-2} & -1 & 6 & -1 \\
			0 & 0 & 0 & 0 & 0 & -z^{-1} - z^{-3} - z^3 - z & 0 & 0 & -1 - z^2 - z^4 - z^{-2} & 8 
		\end{array}
		\right).  
		$$
	\end{example}

\begin{table}[t]
	{\scriptsize
		\begin{center}
			\begin{tabular}{|c|cccccccccc| }
				\hline
				$\zeta=e^{i\frac{2\pi}{8}}$, $z=\zeta^r$ & $\lambda_{r,1}$  & $\lambda_{r,2}$  & $\lambda_{r,3}$  & $\lambda_{r,4}$  & $\lambda_{r,5}$& $\lambda_{r,6}$& $\lambda_{r,7}$& $\lambda_{r,8}$& $\lambda_{r,9}$& $\lambda_{r,10}$ \\
				\hline\hline
				$\spec(\B(\zeta^0))$ & \bf 0 & 1.506 & 3.246 &  4 & 4 &4.890 & 5.452 & 6 & 7.604 & 11.30  \\
				\hline
				$\spec(\B(\zeta^1))=\spec(\B(\zeta^7))$ & \bf 0.586 & 2.215 & 3.126 & $4^*$ &  4.586 & 5.025 & 5.257 & 6.288 & $8^*$ & 8.917  \\
				\hline
				$\spec(\B(\zeta^2))=\spec(\B(\zeta^6))$ & 0.949 & \bf 2  &  2.764 & 3.097 & 4.5173 & 5.194 & 6.534 & 7.230 & 7.709 & $8^*$  \\
				\hline
				$\spec(\B(\zeta^3))=\spec(\B(\zeta^5))$ &  1.108 &  1.712  &  \bf 3.414 &  3.469 &  $4^*$ & 4.874 & 5.718 & 7.414 & $8^*$ & 8.290  \\
				\hline
				$\spec(\B(\zeta^4))$ & 1.079 & 1.330 & 2.0  & \bf 4 & 4 &  4  &  5.522 &   6.403 & 9.34 & 10.257\\
				\hline
			\end{tabular}
	\end{center}}
	\caption{All the eigenvalues of the matrices $\B(\zeta^r)$, which yield the eigenvalues of the 4-token graph $F_4(C_8)$ plus six 8's and four 4's (those marked with `*'). The values in boldface correspond to the eigenvalues of $C_8$.}
	\label{taula:F4C8}
\end{table}

\section{The spectrum of $F_2(C_n)$ through continuous fractions}
\label{sec:cont.frac}

In this section, we concentrate on the case of the $2$-token graph of a cycle and show how to use continuous fractions to derive the spectra of our (Laplacian) polynomial matrices. 
First, using the previous section's method, 
the polynomial matrices of $F_2(C_n)$ can be easily derived. Indeed, for every $z=\zeta^r=e^{ir\frac{2\pi}{n}}$, with $r=0,1,\ldots,n-1$, we get the following $\nu\times \nu$ matrices, $\B_o(z)$ and $\B_e(z)$:

If $n$ is odd, $n=2\nu+1$, 
$$
\B_o(z) =
\left(
\small{
	\begin{array}{cccccc}
		2 & -1-z^{-1} & 0 & 0 & \ldots & 0 \\
		-1-z & 4 & -1-z^{-1} & 0 & \ldots & 0 \\
		0 & -1-z & 4 & -1-z^{-1} & \ddots & 0 \\
		0 & 0 & -1-z & \ddots & \ddots  & 0 \\
		\vdots  & \vdots & \ddots &\ddots & 4 & -1-z^{-1}\\
		0 & 0 & \ldots &0 &-1-z & 4-z^\nu-z^{-\nu}
	\end{array}
}
\right)
$$
and, if $n$ is even, $n=2\nu$,
$$
\B_e(z) =
\left(
\small{
	\begin{array}{cccccc}
		2 & -1-z^{-1} & 0 & 0 & \ldots & 0 \\
		-1-z & 4 & -1-z^{-1} & 0 & \ldots & 0 \\
		0 & -1-z & 4 & -1-z^{-1} & \ddots & 0 \\
		0 & 0 & -1-z & \ddots & \ddots  & 0 \\
		\vdots  & \vdots & \ddots &\ddots & 4 & -1-z^{-1}\\
		0 & 0 & \ldots &0 &-1-z-z^{\nu}-z^{\nu+1} & 4
	\end{array}
}
\right).
$$

The authors also obtained these matrices \cite{dfs19} by using some ad-hoc geometrical and symmetrical reasoning.
Notice that both matrices are tridiagonal. This is because the base graph of $F_2(C_n)$, seen as a lift graph, is path-shaped.
	
	In what follows, every vertex of $F_2(C_n)$ is represented with an ordered pair $(i,j)$, where
	$i,j\in \Z_n$, with $i+h=j$, and $h=i-j=\dist(i,j)>0$ in $C_n$. 
	For instance, the vertices of $F_2(C_6)$ with one $0$
	are $(1,0)$, $(2,0)$, $(3,0)\equiv (0,3)$, $(0,4)$, and $(0,5)$,
	and those of $F_2(C_7)$ are $(1,0)$, $(2,0)$, $(3,0)$, $(0,4)$, $(0,5)$, and $(0,6)$.
	
	The following result, used by the authors in \cite{rdfm23}, is a particular case of Theorem \ref{th:sp-lifts}.
	
	\begin{proposition}[\cite{rdfm23}]
		\label{propo:eigenvec}
		Let $\L$ be the Laplacian matrix of $F_2(C_n)$.
		Every eigenvalue $\lambda$ of $\L$, with $n=2\nu+1$ or $n=2\nu$, has an eigenvector $\x\in \Re^{{n\choose 2}}$ with components
		\begin{equation}
			\label{vector-y-2}
			x_{(i,j)}=f_{i-j}\zeta^j=f_{h}\zeta^j\qquad j=0,\ldots,n-1,\ h=1,\ldots,\nu,
		\end{equation}
		where $\zeta$ is an $n$-th root of unity, and $\f= (f_1,\ldots,f_{\nu})$ is a $\lambda$-eigenvector of the matrix $\B_o(\zeta)$ if $n$ is odd, and $\B_e(\zeta)$ if $n$ is even. Moreover, in the latter case, if $r\neq n/2$, we have that
		$f_i\neq 0$ for $i=1,\ldots,\nu-1$.
	\end{proposition}
	
	In the following result, we use continuous fractions to determine the whole spectrum of $F_2(C_n)$ from the above polynomial matrices.
	
	\begin{theorem}
		\label{th:basic}
		For a given odd $n=2\nu+1$ or even $n=2\nu$, let $Z=\frac{4-\lambda}{2\cos(r\pi/n)}$, and $\alpha=\frac{1}{\cos(r\pi/n)}$, where $r\neq n/2$ in the even case.
		From some  $Q_{\nu-1}=\Psi(Z)$,
		consider the continuous fraction $Q_1=Q_1(\lambda)$ defined recursively as
		\begin{equation}
			Q_{\nu-2}=\frac{1}{Z-Q_{\nu-1}},\ 
			\  \ldots\ ,Q_2=\frac{1}{Z-Q_3},\ Q_1=\frac{1}{Z-Q_2}.
			\label{recur-odd}
		\end{equation}
		That is,
		$$ 
		Q_1=\frac{1}{Z-\frac{1}{Z-\frac{1}{\ddots\ \frac{1}{ Z-\Psi(Z)}}}}.
		$$
		Then, 
		\begin{itemize}
			\item[$(a)$]
			If $n$ is odd, $n=2\nu+1$, and $Q_{\nu-1}=\frac{1}{Z-(-1)^r}$, the values of $\lambda$ that are the solutions of the equation $Q_1=Z-\alpha$, for $r\in\{0,1,\ldots,n-1\}$, are all the  ${n\choose 2}$ eigenvalues of $F_2(C_n)$.
			\item[$(b.1)$]
			If $n$ is even, $n=2\nu$, $r$ is even, $r\neq n/2$, and $Q_{\nu-1}=\frac{2}{Z}$, the solutions of the equation $Q_1=Z-\alpha$, for $r\in\{0,2,\ldots,n-2\}\backslash \{n/2\}$,  are $\nu(\frac{n}{2}-1)$ eigenvalues of $F_2(C_n)$.
			\item[$(b.2)$]
			If $n$ is even and $r=n/2$, then $\B(z)$ is a diagonal matrix with entries $2,4,\stackrel{(\nu-1)}{\ldots},4$, which correspond to $\nu$ eigenvalues of  $F_2(C_n)$.
			\item[$(b.3)$]
			If $n$ is even, $n=2\nu$, $r$ is odd, and $Q_{\nu-2}=\frac{1}{Z}$ (that is, $Q_{\nu-1}=0$), the solutions of the equation $Q_1=Z-\alpha$, for $r\in\{1,3,\ldots,n-1\}$, are $(\nu-1)\frac{n}{2}$ eigenvalues of $F_2(C_n)$ different from $\lambda=4$.
		\end{itemize}
	\end{theorem}
	
	\begin{proof}
		Let $\L$ be the Laplacian matrix of $F_2(C_n)$, and $R(n)=\{e^{ir\frac{2\pi}{n}}:r=0,\ldots,n-1\}$ the set of $n$-th roots of unity. 
		
		$(a)$ From Proposition \ref{propo:eigenvec}, a $\lambda$-eigenvector $\x$ of $\L$
		has entries
		\begin{equation}
			\x_{(i,j)}=f_h \zeta^j,\quad \zeta\in R(n), \ h=1,\ldots,\nu,
		\end{equation}
		where $\f=(f_1,\ldots,f_{\nu})$ is a $\lambda$-eigenvector of $\B_o(\zeta)$.
		Assuming that $f_1,\ldots,f_{\nu-1}$ are different from zero, let $Q_h(r)=\frac{f_{h+1}}{f_h}\zeta^{-1/2}=\frac{f_{h+1}}{f_h}e^{-ir\frac{\pi}{n}}$ for $h=1,\ldots,\nu-1$. We now distinguish three cases:
		
		\begin{itemize}
			\item[$(i)$] $h=1$.
			The equality $\B_o(\zeta)\f=\lambda \f$ gives
			$(2-\lambda)f_1=f_2(1+\zeta^{-1})$. Thus, by dividing by $f_1$ and taking common factor $\zeta^{-1/2}$, we get $(2-\lambda)=\frac{f_2}{f_1}\zeta^{-1/2}(\zeta^{1/2}+\zeta^{-1/2})$. That is,
			\begin{equation}
				\label{Q1}
				Q_1=\frac{2-\lambda}{\zeta^{1/2}+\zeta^{-1/2}}=\frac{2-\lambda}{2\cos(r\pi/n)}=Z-\alpha.
			\end{equation}
			\item[$(ii)$] $1<h<\nu$.
			The equality $\B_o(\zeta)\f=\lambda \f$ gives
			$(4-\lambda)f_h=f_{h-1}(1+\zeta)+f_{h+1}(1+\zeta^{-1})$. Dividing now by $f_h$ and taking common factor $\zeta^{-1/2}$, we get 
			$4-\lambda=\frac{f_{h-1}}{f_h}\zeta^{1/2}(\zeta^{-1/2}+\zeta^{1/2})+
			\frac{f_{h+1}}{f_h}\zeta^{-1/2}(\zeta^{1/2}+\zeta^{-1/2})$. That is,
			\begin{equation}
				\label{Qh-Qh-1}
				\frac{4-\lambda}{\zeta^{1/2}+\zeta^{-1/2}}=\frac{4-\lambda}{2\cos(r\pi/n)}=Z=\frac{1}{Q_{h-1}}+Q_{h}
			\end{equation}
			and, hence,
			$$
			Q_{h-1}=\frac{1}{Z-Q_h}\quad \mbox{for $h=2,\ldots,\nu-1$.}
			$$
			\item[$(iii)$] $h=\nu=(n-1)/2$.
			Now, the last equality of $\B_o(\zeta)\f=\lambda \f$ leads to
			$(4-\lambda)f_{\nu}=f_{\nu-1}(1+\zeta)+f_{\nu}(\zeta^{\nu}+\zeta^{-\nu})$. Dividing  by $f_{\nu}$ and taking common factor $\zeta^{1/2}$, we get 
			$4-\lambda=\frac{f_{\nu-1}}{f_\nu}\zeta^{1/2}(\zeta^{-1/2}+\zeta^{1/2})+
			(\zeta^{\nu}+\zeta^{-\nu})$. That is,
			\begin{align}
				Q_{\nu-1}&=\frac{\zeta^{1/2}+\zeta^{-1/2}}
				{4-\lambda-(\zeta^{\nu}+\zeta^{-\nu})}=\frac{2\cos(r\pi/n)}{4-\lambda-2\cos(r\nu 2\pi/n)} \nonumber\\
				&=\frac{2\cos(r\pi/n)}{4-\lambda-(-1)^r 2\cos(r\pi/n)}=\frac{1}{Z-(-1)^r},\label{Qnu}
			\end{align}
			since $\cos(r\nu 2\pi/n)=\cos\left(\frac{r\pi(n-1)}{n}\right)=\cos(r\pi)\cos(r\pi/n)$ and $4-\lambda=2Z\cos(r\pi/n)$.
		\end{itemize}
		
		$(b.1-3)$
		Let $n=2\nu$.
		Using the same notation as in  $(i)$, the cases $h=1$ and $1<\nu$ lead to the same results as before, provided that $r\neq n/2$. That is, $Q_1=\frac{2-\lambda}{2\cos(r\pi/n)}=Z-\alpha$, and $Q_{h-1}=\frac{1}{Z-Q_h}$ for $h=2,\ldots,\nu-1$.
		
		In contrast, for $h=\nu=n/2$, the last equality of $\B_e(\zeta)\f=\lambda\f$ gives
		\begin{align}
			(4-\lambda)f_{\nu}& =f_{\nu-1}(1+\zeta+\zeta^{\nu}+\zeta^{\nu+1})
			=f_{\nu-1}\zeta^{1/2}(\zeta^{-1/2}+\zeta^{1/2}+\zeta^{\nu-1/2}+\zeta^{\nu+1/2})\nonumber\\
			&=f_{\nu-1}\zeta^{1/2}\left[2\cos\left(\frac{r\pi}{n}\right)+2\cos\left(\frac{r(n-1)\pi}{n}\right)\right]
			\nonumber\\
			&=f_{\nu-1}\zeta^{1/2}2\cos\left(r\pi/n\right)[1+(-1)^r].\label{4-Lfnu}
		\end{align}
		Thus, dividing by $f_{\nu}$, and solving for $Q_{\nu-1}$, we have
		\begin{align}
			\label{Qnu-1}
			Q_{\nu-1} &=\frac{2\cos\left(r\pi/n\right)[1+(-1)^r]}{4-\lambda}.
		\end{align}
		Finally, we must distinguish three cases:
		\begin{enumerate}
			\item[$(b.1)$]
			If $r$ is even, with $r\neq n/2$, $Z$ and $\alpha$ are well defined, and \eqref{Qnu-1} gives
			$$ 
			Q_{\nu-1} =\frac{4\cos\left(r\pi/n\right)}{4-\lambda}=\frac{2}{Z},
			$$
			as claimed.
			\item[$(b.2)$] 
			If $r=n/2$, $\zeta=\zeta^{-1}=\cos(\pi)=-1$ and, hence $\B_e(\zeta)=\diag(2,4,\ldots,4)$ yields the mentioned eigenvalues of $F_2(C_n)$.
			\item[$(b.3)$]
			If $r$ is odd, we must go back to \eqref{4-Lfnu}, giving $(4-\lambda)f_{\nu}=0$. Thus, either $\lambda=4$ or 
			$f_{\nu}=0$. In the first case, such a value does not appear as an eigenvalue of $F_2(C_n)$, for the same reason explained in Example \ref{ex:F3(C6)}. Namely, the last component of the eigenvector of $\B_o(z)$ is not zero, and $o(r)\nmid 2$. In the second case, we need to look at the $(\nu-1)$-th equality of $\B_e(\zeta)\f=\lambda\f$, which yields
			$$
			(4-\lambda)f_{\nu -1}=f_{\nu-2}(1+\zeta)=f_{\nu-2}\zeta^{1/2}2\cos\left(r\pi/n\right).
			$$
			Hence, 
			$$
			Q_{\nu-2}=\frac{2\cos\left(r\pi/n\right)}{4-\lambda}=\frac{1}{Z},
			$$
			as stated. More directly, in this case, \eqref{Qnu-1} gives $Q_{\nu-1}(4-\lambda)=0$ and, hence, if $\lambda\neq 4$, we have $Q_{\nu-1}=0$ and $Q_{\nu-2}=\frac{1}{Z}$.
		\end{enumerate}
		Notice that the number of eigenvalues obtained in the cases $(b.1)-(b.3)$ is $\nu(\frac{n}{2}-1)+\nu+(\nu-1)\frac{n}{2}$ which, with $\nu=n/2$, gives a total of ${n\choose 2}$, as expected.
		This completes the proof.
	\end{proof}
	
	\subsection{Some examples}
	\begin{itemize}
		\item 
		For $n=7$, with $F_2(C_7)$ shown on the left-side of Figure \ref{F2(C7)+F3(C7)}, the rational function $Q_1$ turns out to be
		$$
		Q_1(Z)=\frac{Z-(-1)^r}{Z^2-(-1)^r Z-1}.
		$$
		Solving for $\lambda$ the equation $Q_1(Z)=Z-\alpha$ with $\alpha= \frac{1}{\cos(r\pi/7)}$ and $Z=\frac{4-\lambda}{2\cos(r\pi/7)}$, with $r=0,\ldots,6$  we obtain all the eigenvalues of $F_2(C_7)$, see Table \ref{taula:C7}.
		
		\begin{table}[t]
			\begin{center}
				\begin{tabular}{|c|ccc| }
					\hline
					$\zeta=e^{i\frac{2\pi}{7}}$, $z=\zeta^r$ & $\lambda_{r,1}$  & $\lambda_{r,2}$  & $\lambda_{r,3}$   \\
					\hline\hline
					$\spec(\B(\zeta^0))$ & \bf 0 & 2.0 &  6.0   \\
					\hline
					$\spec(\B(\zeta^1))=\spec(\B(\zeta^6))$ & \bf 0.7530   &  3.9363  &  7.1125 \\
					\hline
					$\spec(\B(\zeta^2))=\spec(\B(\zeta^5))$ &  1.1633 & \bf 2.4450 & 5.1446   \\
					\hline
					$\spec(\B(\zeta^3))=\spec(\B(\zeta^4))$ & 1.9019 & \bf 3.8019  &  4.7411   \\
					\hline
				\end{tabular}
			\end{center}
			\caption{All the eigenvalues of the matrices $\B(\zeta^r)$, which yield the eigenvalues of the 2-token graph $F_2(C_7)$. The values in boldface correspond to the eigenvalues of $C_7$.}
			\label{taula:C7}
		\end{table}
		
		\item
		For $n=8$, the rational functions $Q_1$ turn out to be
		$$
		Q_1(Z)=\frac{Z^2-2}{Z(Z^2-3)}\quad \mbox{for $r$ even, $r\neq 4$, and}\quad Q_1(Z)=\frac{Z}{Z^2-1}\quad\mbox{for $r$ odd.}
		$$
		In this case, solving for $\lambda$ the equations in Theorem \ref{th:basic}$(b)$ (or looking at the diagonal matrix $\B_e$ when $r=n/2$), we get all the eigenvalues of $F_2(C_8)$, see Table \ref{taula:C8}.
	\end{itemize}
	
	\begin{table}[t!]
		\begin{center}
			\begin{tabular}{|c|cccc|}
				\hline
				$\zeta=e^{i\frac{2\pi}{8}}$, $z=\zeta^r$ & $\lambda_{r,1}$  & $\lambda_{r,2}$  & $\lambda_{r,3}$  & $\lambda_{r,4}$   \\
				\hline\hline
				$\spec(\B(\zeta^0))$ & \bf 0 & 1.5060 & 4.8900 &  7.60387   \\
				\hline
				$\spec(\B(\zeta^1))=\spec(\B(\zeta^7))$ & \bf 0.5857 & 3.1259 & $4.0^*$ &  6.2882  \\
				\hline
				$\spec(\B(\zeta^2))=\spec(\B(\zeta^6))$ & 0.9486 & \bf 2.0  &  4.5173 & 6.5340  \\
				\hline
				$\spec(\B(\zeta^3))=\spec(\B(\zeta^5))$ &  1.7117  &  \bf 3.4142 &  $4.0^*$  & 4.8740  \\
				\hline
				$\spec(\B(\zeta^4))$ & 2.0   & \bf 4.0  &  4.0 &   4.0 \\
				\hline
			\end{tabular}
		\end{center}
		\caption{All the eigenvalues of the matrices $\B(\zeta^r)$, which yield the eigenvalues of the 2-token graph $F_2(C_8)$ plus four 4's (those marked with `*'). The values in boldface correspond to the eigenvalues of $C_8$.}
		\label{taula:C8}
	\end{table}
	
	\section{The characteristic polynomials of $F_2(C_n)$}
	\label{sec:charac-pol}
	The results of the previous section allow us to give a closed formula for the characteristic polynomial of $F_2(C_n)$.
	Let us first consider the case of odd $n$.
	
	\begin{theorem}
		\label{th:n-odd}
		Given $n=2\nu+1\ge 5$, let $\L$ be the Laplacian matrix of $F_2(C_n)$. Let $Z=\frac{4-\lambda}{2\cos(r\pi/n)}$, $\alpha=\frac{2}{2\cos(r\pi/n)}$, and $\rho_{1,2}=\frac{1}{2}(Z\pm\sqrt{Z^2-4})$.  Then, up to a multiplicative constant, the characteristic polynomial of $\L$  is
		\begin{equation}
			\phi_{\small \L}(\lambda)=\prod_{r=0}^{n-1}\frac{1}{\sqrt{Z^2-4}}\left[(\rho_2-(-1)^r)(\rho_2-\alpha)\rho_2^{\nu-1}
			-(\rho_1-(-1)^r)(\rho_1-\alpha)\rho_1^{\nu-1}\right].
		\end{equation}
	\end{theorem}
	
	\begin{proof}
		Let $Q_i=\frac{R_{\nu-i-1}}{S_{\nu-i-1}}$ for $i=1,\ldots,\nu-1$. Then, starting with $\frac{R_0}{S_0}=\frac{1}{Z-(-1)^r}$, the recurrence in  \eqref{recur-odd} is
		$$
		\frac{R_i}{S_i}=
		\frac{1}{Z-\frac{R_{i-1}}{S_{i-1}}}
		=\frac{S_{i-1}}{ZS_{i-1}-R_{i-1}}.
		$$
		This holds when $R_i=S_{i-1}$ and $S_i=ZS_{i-1}-R_{i-1}$
		or, in matrix form,
		$$
		\left(
		\begin{array}{c}
			R_i \\
			S_i
		\end{array}
		\right)=
		\left(
		\begin{array}{cc}
			0 & 1\\
			-1   & Z
		\end{array}
		\right)
		\left(
		\begin{array}{c}
			R_{i-1} \\
			S_{i-1}
		\end{array}
		\right)=\cdots=\M^i\left(
		\begin{array}{c}
			R_{0} \\
			S_{0}
		\end{array}
		\right)=
		\M^i\left(
		\begin{array}{c}
			1 \\
			Z-(-1)^r
		\end{array}
		\right),
		$$
		where $\M=\left(
		\begin{array}{cc}
			0 & 1\\
			-1   & Z
		\end{array}
		\right)$. Therefore, the equality in which the solutions are the eigenvalues of $\B(\zeta)$ (with $\zeta=e^{ir\frac{2\pi}{n}}$), namely, $Q_1=\frac{R_{\nu-2}}{S_{\nu-2}}-(Z-\alpha)=0$, becomes
		\begin{equation}
			0=R_{\nu-2}-(Z-\alpha)S_{\nu-2} 
			=\left(
			1,  -Z+\alpha
			\right)\M^{\nu-2}\left(
			\begin{array}{c}
				1 \\
				Z-(-1)^{r}
			\end{array}
			\right).
			\label{eq:charac}
		\end{equation}
		
		Diagonalizing $\M$, which has eigenvalues
		$$
		\rho_1=\frac{1}{2}(Z+\sqrt{Z^2-4})\quad\mbox{and}\quad \rho_2=\frac{1}{2}(Z-\sqrt{Z^2-4}),
		$$
		we get 
		$$
		\U^{-1}\M\U=\D=\left(
		\begin{array}{cc}
			\rho_1 & 0\\
			0 & \rho_2
		\end{array}
		\right)
		$$
		with
		$$
		\U=\left(
		\begin{array}{cc}
			1 & 1\\
			\rho_1 & \rho_2
		\end{array}
		\right)
		\quad \mbox{ and }\quad
		\U^{-1}=\frac{1}{\rho_1-\rho_2}\left(
		\begin{array}{rr}
			-\rho_2 & 1\\
			\rho_1 &  -1
		\end{array}
		\right).
		$$
		Thus, plugging $\M=\U\D\U^{-1}$ in \eqref{eq:charac}, we get the equation
		\begin{align}
			& \left(
			1,  -Z+\alpha
			\right)\U\D^{\nu-2}\U^{-1}\left(
			\begin{array}{c}
				1 \\
				Z-(-1)^{r}
			\end{array}
			\right) \label{pol-charac-n-odd}\\
			&=\left(
			1,  -Z+\alpha
			\right)\left(
			\begin{array}{cc}
				1 & 1\\
				\rho_1 & \rho_2
			\end{array}
			\right)\left(
			\begin{array}{cc}
				\rho_1^{\nu-2} & 0\\
				0 & \rho_2^{\nu-2}
			\end{array}
			\right)\frac{1}{\rho_1-\rho_2}\left(
			\begin{array}{rr}
				-\rho_2 & 1\\
				\rho_1 &  -1
			\end{array}
			\right)\left(
			\begin{array}{c}
				1 \\
				Z-(-1)^{r}
			\end{array}
			\right)\nonumber\\
			&=\frac{1}{\sqrt{Z^2-4}}\left[(\rho_2-(-1)^r)(\rho_2-\alpha)\rho_2^{\nu-1}
			-(\rho_1-(-1)^r)(\rho_1-\alpha)\rho_1^{\nu-1}
			\right]=0,
			\label{eq:charac2}
		\end{align}
		where we used that $\rho_1-\rho_2=\sqrt{Z^2-4}$, $\rho_1\rho_2=1$, $Z-\rho_1=\rho_2$, and $Z-\rho_2=\rho_1$.
		
		Then, the non-zero term in \eqref{eq:charac2} is, up to a multiplicative constant, the characteristic polynomial of $\B(r)=\B(\zeta)$ with $\zeta=e^{ir\frac{2\pi}{9}}$, and the result follows since, by Theorem \ref{th:basic}$(a)$, $\spec \L=\bigcup_{r=0}^{n-1} \spec\B(r)$.
	\end{proof}
	
\begin{figure}[t]
		\begin{center}
			\includegraphics[width=9cm]{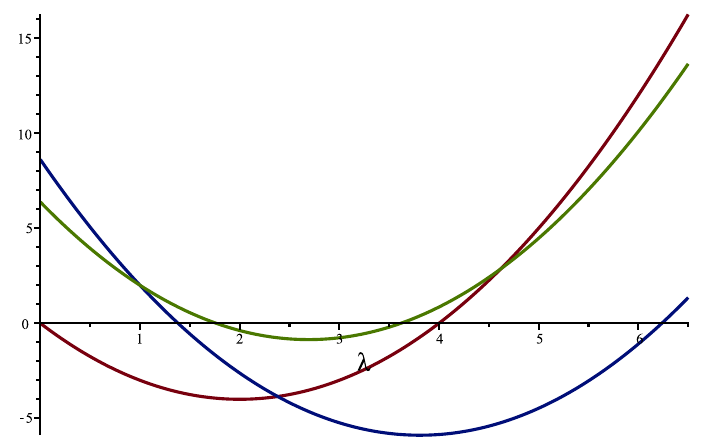}
			\caption{The characteristic polynomials of $\B(r)$, with $r=0,1,2$, for $F_2(C_5)$.}
			\label{charac-pols-F_2(C_5)}
		\end{center}
	\end{figure}

\begin{figure}[t]
		\begin{center}
			\includegraphics[width=9cm]{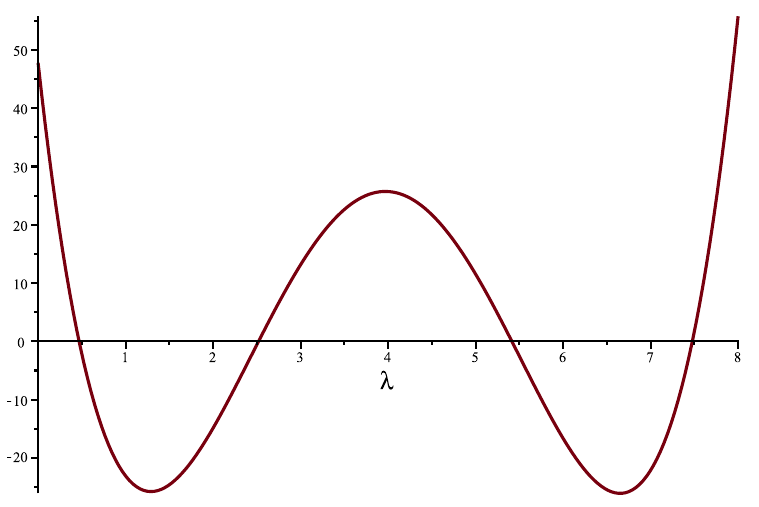}
			\caption{The characteristic polynomial of $\B(1)$, with smallest root being the algebraic connectivity of $F_2(C_9)$. }
			\label{charac-pol-F_2(C_9)}
		\end{center}
	\end{figure}

	Let us see two examples.
	\begin{itemize}
		\item
		When $n=5$, the characteristic polynomials $\phi_r$ of $\B(r)$, for $r=0,\ldots,4$ are the following:
		\begin{align*}
			\phi_{0}(\lambda) &= \lambda^2-4\lambda,\\
			\phi_{1}(\lambda)=\phi_4(\lambda) &= \lambda^2 + (-\sqrt{5}/2 - 13/2)\lambda + 15/2 + \sqrt{5}/2,\\
			\phi_{2}(\lambda)=\phi_3(\lambda) &=\lambda^2 + (\sqrt{5}/2 - 13/2)\lambda + 15/2 - \sqrt{5}/2.
		\end{align*}
		See Figure \ref{charac-pols-F_2(C_5)}.
		The smallest root of $\B(1)=\B(4)$, is $\frac{1}{2}(5-\sqrt{5})\approx 1.3819$, corresponding to the algebraic connectivity of $F_2(C_5)$.
		\item
		When $n=9$ and $r=1$, the characteristic polynomial $\phi_1$ of $\B(1)$  (coefficients with four digits of approximation) is:
		$$
		\phi_1(\lambda)=
		\lambda^4 - 15.88\lambda^3 + 80.19\lambda^2 - 136.2\lambda + 47.79,
		$$
		with smallest root being $\alpha(F_2(C_9))\approx 0.4679$, see Figure 
		\ref{charac-pol-F_2(C_9)}.
	\end{itemize}
	
	Now, we consider the case of $n$ even.
	\begin{theorem}
		Given $n=2\nu\ge 4$, let $\L$ be the Laplacian matrix of $F_2(C_n)$. Let $Z=\frac{4-\lambda}{2\cos(r\pi/n)}$, $\alpha=\frac{2}{2\cos(r\pi/n)}$, and $\rho_{1,2}=\frac{1}{2}(Z\pm\sqrt{Z^2-4})$.  Then, up to a multiplicative constant, the characteristic polynomial of $\L$ is
		\begin{align}
			\phi_{\small \L}(\lambda)& =\prod_{\stackrel{r=0}{r\ even, r\neq n/2}}^{n-2}
			(1-(Z-\alpha)\rho_2)\rho_2^{\nu-2}
			+(1-(Z-\alpha)\rho_1)\rho_1^{\nu-2}\nonumber\\
			&\times \prod_{\stackrel{r=1}{r\ odd}}^{n-1}
			\frac{1}{\sqrt{Z^2-4}}\left[(-1+(Z-\alpha)\rho_2)\rho_2^{\nu-2}
			+(1+(Z-\alpha)\rho_1)\rho_1^{\nu-2}\right]\nonumber\\
			&\times (\lambda-2)(\lambda-4)^{\nu-1}.
		\end{align}
	\end{theorem}
	
	\begin{proof}
		The proof goes along the same lines of reasoning as in Theorem \ref{th:n-odd}. The only changes are the following:
		\begin{itemize}
			\item 
			When $r$ is even, $r\neq n/2$, the expression in \eqref{pol-charac-n-odd} becomes
			$$
			\left(
			1,  -Z+\alpha
			\right)\U\D^{\nu-2}\U^{-1}\left(
			\begin{array}{c}
				2 \\
				Z
			\end{array}
			\right).
			$$
			\item 
			When $r$ is odd, the expression in \eqref{pol-charac-n-odd} becomes
			$$
			\left(
			1,  -Z+\alpha
			\right)\U\D^{\nu-3}\U^{-1}\left(
			\begin{array}{c}
				1 \\
				Z
			\end{array}
			\right).
			$$
			\item 
			When $r=n/2$, $\B(r)=\diag(2,4,\ldots,4)$ and, hence \eqref{pol-charac-n-odd} must be $(\lambda-2)(\lambda-4)^{\nu-1}$.
		\end{itemize}
	\end{proof}
	
\begin{figure}[t]
		\begin{center}
			\includegraphics[width=9cm]{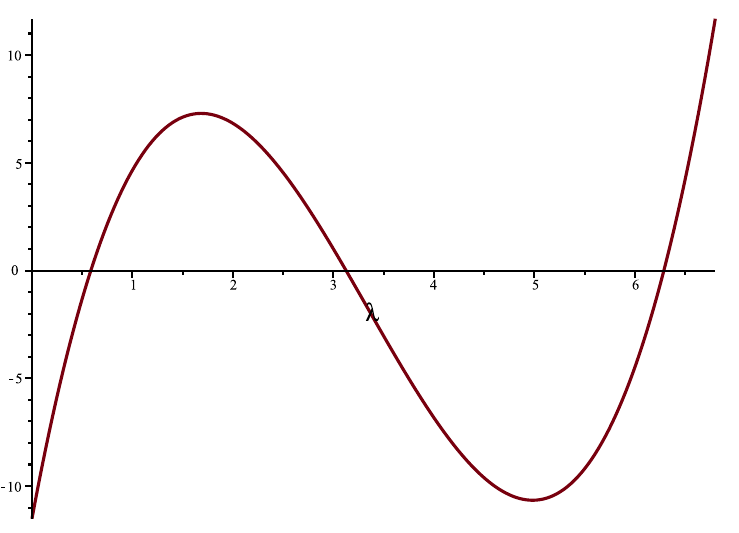}
			\caption{The characteristic polynomial of $\B(1)$, with the smallest root being the algebraic connectivity of $F_2(C_8)$. }
			\label{charac-pol-F_2(C_8)}
		\end{center}
	\end{figure}

	For instance, when $n = 8$ and $r = 1$, the characteristic polynomial of $\B(1)$ (coefficients
	with two digits of approximation) is:
	$$
	\phi_{1}(\lambda) =\lambda^3 - 10\lambda^2 + 25.17\lambda - 11.51,
	$$
	with the smallest root being $\alpha(F_2(C_8)) \approx 0.5857$, see Figure \ref{charac-pol-F_2(C_8)}.
	
\end{document}